\documentclass[10pt,reqno]{amsart}
\usepackage{amsmath,amssymb,amsthm}
\usepackage{graphicx}
\usepackage[all]{xypic}
\usepackage[utf8,nocaptions]{vietnam}
\usepackage{hyperref}
\input xypic
\def\thmsection{section}
\def\thmchangesection{changesection}

\def\thmchangechapter{changechapter}
\def\thmchange{change}
\def\thmplain{plain}
\ifx\theoremnumberstyle\thmplain
  \theoremstyle{break-italic}
  \newtheorem{satz}{Satz}
\else
  \ifx\theoremnumberstyle\thmsection
    \theoremstyle{break-italic}
    \newtheorem{satz}{Satz}[section]
  \else
    \ifx\theoremnumberstyle\thmchange
      \swapnumbers
      \theoremstyle{break-italic}
      \newtheorem{satz}{Satz}
    \else
      \ifx\theoremnumberstyle\thmchangesection
         \swapnumbers
         \theoremstyle{break-italic}
         \newtheorem{satz}{Satz}[section]
      \else
        \ifx\theoremnumberstyle\thmchangechapter
           \swapnumbers
           \theoremstyle{break-italic}
           \newtheorem{satz}{Satz}[chapter]
        \else
          \ifx\theoremnumberstyle\thmsection
             \theoremstyle{break-italic}
             \newtheorem{satz}{Satz}[section]
          \else
            \theoremstyle{break-italic}
            \newtheorem{satz}{Satz}[section]
          \fi
        \fi
      \fi
    \fi
  \fi
\fi

\theoremstyle{break-italic}
\newtheorem{theorem}[satz]{Theorem}
\newtheorem{lemma}[satz]{Lemma}

\newtheorem{corollary}[satz]{Corollary}
\newtheorem{Proposition}[satz]{Proposition}

\newtheorem*{conjecture*}{Conjecture}

\theoremstyle{break-roman}
\newtheorem{definition}[satz]{Definition}

\newtheorem{example}[satz]{Example}

\newtheorem{remark}[satz]{Remark}

\theoremstyle{standard}

\newtheorem*{claim}{Claim}

\theoremstyle{varthm-roman}
\newtheorem*{varthm-roman}{}
\theoremstyle{varthm-italic}
\newtheorem*{varthm-italic}{}
\theoremstyle{varthm-roman-break}
\newtheorem*{varthm-roman-break}{}
\theoremstyle{varthm-italic-break}
\newtheorem*{varthm-italic-break}{}
\theoremstyle{varthm-roman-no-punctuation}
\newtheorem{varthm-roman-no-punctuation-numbered}[satz]{}
\theoremstyle{varthm-italic-no-punctuation}
\newtheorem{varthm-italic-no-punctuation-numbered}[satz]{}

\newenvironment{varthm-roman-numbered}[1]{
  \begin{varthm-roman-no-punctuation-numbered}
    \mbox{\rm\textbf{#1}}
  }{\end{varthm-roman-no-punctuation-numbered}}
\newenvironment{varthm-italic-numbered}[1]{
  \begin{varthm-italic-no-punctuation-numbered}
    \mbox{\rm\textbf{#1}}
  }{\end{varthm-italic-no-punctuation-numbered}}
\newenvironment{varthm-roman-break-numbered}[1]{
  \begin{varthm-roman-no-punctuation-numbered}
    \mbox{\rm\textbf{#1}\newline}
  }{\end{varthm-roman-no-punctuation-numbered}}
\newenvironment{varthm-italic-break-numbered}[1]{
  \begin{varthm-italic-no-punctuation-numbered}
    \mbox{\rm\textbf{#1}}\newline
  }{\end{varthm-italic-no-punctuation-numbered}}
\numberwithin{equation}{section}

\def\ex{\begin{example}
  }
  \def\eex{\end{example}}
\def\thr{\begin{theorem}}
\def\ethr{\end{theorem}}
\def\pro{\begin{Proposition}}
\def\epro{\end{Proposition}}
\def\coro{\begin{corollary}}
\def\ecoro{\end{corollary}}
\def\df{\begin{definition}}
\def\edf{\end{definition}}
\def\lm{\begin{lemma}}
\def\elm{\end{lemma}}
\def\pf{\begin{proof}}
\def\epf{\end{proof}}
\def\problem{\begin{problem}}
\def\eproblem{\end{problem}}
\def\dlim{\displaystyle\lim}

\def\it{\begin{itemize}}
\def\hit{\end{itemize}}
\def\rem{\begin{remark}}
\def\erem{\end{remark}}

\def\cla{\begin{claim}}
\def\ecla{\end{claim}}
\def\dlim{\displaystyle\lim}

\renewcommand{\tt}{\exists}
\newcommand{\vc}{\infty}
\newcommand{\vm}{\forall}
\newcommand{\tru}{\setminus}

\newcommand{\Sr}{\Longrightarrow}

\newcommand{\ps}{\dfrac}
\newcommand{\can}{\sqrt}
\newcommand{\mtn}{\rightarrow}

\newcommand{\thuc}{\Bbb  R}

\def\n{\Bbb N}

\def\i{\mbox{\bf I}}
\def\hoa{\mathcal}

\begin{document}
\title[Some Positivstellens\"{a}tze for polynomial matrices] {Some Positivstellens\"{a}tze for polynomial matrices}

\author{L\^{e} C\^{o}ng-Tr\`{i}nh }
\address{Department of Mathematics, Quy Nhon University\\
170 An Duong Vuong, Quy Nhon, Binh Dinh, Vietnam}

\email{lecongtrinh@qnu.edu.vn}


\begin{abstract} In this paper we give a version of Krivine-Stengle's Positivstellensatz, Schweighofer's Positivstellensatz, Scheiderer's local-global principle, Scheiderer's Hessian criterion and  Marshall's boundary Hessian conditions   for polynomial matrices, i.e. matrices with entries from the ring of polynomials in the variables \(x_1,\cdots,x_d\) with real coefficients. Moreover, we characterize 
Archimedean quadratic modules of  polynomial matrices, and  study the relationship between the compactness of a subset   in \(\mathbb R^{d}\) with respect to a subset $\hoa{G}$ of polynomial matrices  and the Archimedean property of the preordering and the quadratic module generated by $\hoa{G}$.
\end{abstract}
\maketitle

\section{Introduction}

Let $\thuc[X]:=\thuc[x_1,\cdots,x_d]$ be the ring of polynomials in the variables $x_1,\cdots,x_d$ with real coefficients. Denote by $\sum \thuc[X]^{2}$ the set of \textit{sums of squares}  in $\thuc[X]$, i.e. the set of finite sums $\sum f_i^{2}, f_i\in \thuc[X]$.   For a subset  $G=\{g_1,\cdots,g_m\}\subseteq \thuc[X]$,  let us consider  the basic closed semi-algebraic set associated to $G$, 
$$ K_G:=\{x\in \thuc^{d} | g_i(x) \geq 0, i=1,\cdots,m\}, $$
the quadratic module generated by $G$,
$$M_G  = \{ t_0+\sum_{i=1}^{m} t_i g_i | t_i \in \sum \thuc[X]^{2}, i=0,1,\cdots,m\}$$
and the preordering generated by $G$,
$$ T_G = \{\sum_{\sigma=(\sigma_1,\cdots,\sigma_m)\in \{0,1\}^{m}} t_{\sigma}g_1^{\sigma_1}\cdots g_m^{\sigma_m}| t_{\sigma}\in \sum \thuc[X]^{2}\}. $$

For a polynomial $f\in \thuc[X]$, it is obvious that if  $f\in M_G$ or $f\in T_G$ then $f(x)\geq 0$ for all $x\in K_G$ (in this case we say $f\geq 0$ on $K_G$).  The converse is in general not true. The Positivstellensatz of Krivine-Stengle (\cite[1964]{Kr}, \cite[1974]{St})  characterizes  polynomials which are positive (resp. non-negative, vanished) on a basic closed semi-algebraic set, but with a "denominator" (for example, $f>0$ on $K_G$ if and only if $pf = 1 + q$ for some $p,q\in T_G$, that is, $f\in \ps{1}{p}(1+ T_G)$ with denominator $p\in T_G$).  \\
A "denominator-free" version of this result is due to Schm\"{u}dgen (1991) which asserts that any  positive  polynomial  on a compact set   $K_G$ belongs to $T_G$. To ensure for $f>0$ on $K_G$ to be in $M_G$, Putinar (1993) required the Archimedean property of $M_G$. Note that the compactness of $K_G$ is equivalent to the Archimedean property of $T_G$ (cf. \cite[Theorem 6.1.1]{Ma}), and if $M_G$ is Archimedean then so is $T_G$,  hence $K_G$ is compact. However the converse is not true in general (see, for example, \cite[Putinar's question, chapter 7]{Ma}).    \\
If $K_G$ is not assumed to be compact, Schweighofer (\cite{Schw}) has given a Positivstellensatz which asserts that if $f\in \thuc[X]$ is a bounded, positive polynomial on $K_G$  and if it has only finitely many asymptotic values on $K_G$ such that all of them are positive then $f\in T_G$.   \\
The case where $K_G$ is compact (resp. $M_G$ is Archimedean),  but $f$ is assumed to have  finitely many  zeros in $K_G$,  Scheiderer (\cite{Sch1}, \cite{Sch2}) has given a Hessian criterion at each zero  of $f$ in $K_G$ for $f$ to be in $T_G$ (resp. $M_G$), using his local-global criterion.   Marshall (cf. \cite{Ma}) has also given  boundary Hessian  conditions at each zero of $f$ in $K_G$ to ensure for $f$ to be in $M_G$. 

The aim of this paper is to study all of these Positivstellens\"{a}tze for polynomial matrices, that is for matrices with entries from $\thuc[X]$.  A matrix version of   Krivine-Stengle's Positivstellensatz was given  by Schm\"{u}dgen (\cite[2009]{Schm}, for non-negative polynomial matrices) and Cimprič (\cite[2012]{Ci}).  Hol-Scherer (\cite[2006]{SchH}, or \cite[2010]{KSch})  has given a matrix version of Putinar's Positivstellensatz. 
Cimprič has also given a version of Schm\"{u}dgen's Positivstellensatz for polynomial matrices in \cite[2013]{CZ}.  

In section 2 we recall  definition of quadratic modules and preorderings in the algebra $\hoa{M}_n(\thuc[X])$ of polynomial matrices, which is proposed by Schm\"{u}dgen (\cite{Schm-1}, \cite{Schm0}, \cite{Schm}) and Cimprič (\cite{Ci0}, \cite{Ci}), and some basic facts used in the paper. In particular, we recall a basic result of Cimprič (see Lemma \ref{lm1.3}) which tells us that any  subset $\bold{K}_\hoa{G}$ of $\thuc^{d}$  associated to $\hoa{G}\subseteq \hoa{S}_n(\thuc[X])$ can be determined again by a subset $G$ of polynomials in $\thuc[X]$ such that the preordering $\hoa{T}_{\hoa{G}}$ (resp. the quadratic module $\hoa{M}_{\hoa{G}}$) contains the preordering $(T_G)^{n}$ (resp. the quadratic module $(M_G)^{n}$). Moreover we recall also  a  basic  result of Schm\"{u}dgen (see Lemma \ref{lm1.4}) which asserts that any symmetric polynomial matrix, together with a square of a non-zero polynomial in $\thuc[X]$, can be diagonalized. This allows us to prove many results of this paper firstly with diagonal matrices, and then with arbitrary symmetric matrices.  

 In section 3 we give a matrix version of Krivine-Stengle's Positivstellensatz (Proposition \ref{pro2.1} and Theorem \ref{thr2.2}). This version for polynomial matrices is simpler than the one given in \cite{Schm} (for positive semidefinite polynomial matrices), however  in general    more complicated than the one given in  \cite{Ci}. But  in our version, the existence of diagonal  polynomial matrices  in the representation of $\bold{F}$ in $(T_G)^{n}$ is  more convenient. 

In section 4 we give a matrix version of Schweighofer's Positivstellensatz (Proposition \ref{pro3.2} and Theorem \ref{thr3.3}). We have a nice representation for  diagonal polynomial matrices, however  in the representation of   an arbitrary symmetric polynomial matrix  we need  a "denominator", namely, a  square of a  non-zero polynomial in $\thuc[X]$ or a conjugation of a matrix in $\hoa{M}_n(\thuc[X])$.

In section 5   we recall definition of  Archimedean quadratic modules in $\hoa{M}_n(\thuc[X])$ and characterize Archimedean  quadratic modules via  the ring of bounded elements with respect to these quadratic modules. We show that the Archimedean property of a quadratic module $M$ in $\thuc[X]$ is the same as that of the quadratic module $M^{n}$ in $\hoa{M}_n(\thuc[X])$, and the compactness of the set $\bold{K}_{\hoa{G}}$ is equivalent to the Archimedean property of the preordering $\hoa{T}_{\hoa{G}}$.  Moreover, we show that if the quadratic module $\hoa{M}_{\hoa{G}}$ of  univariate polynomial matrices is Archimedean then the set $\bold{K}_{\hoa{G}}$ is compact.  

The last section deals with a matrix version of Scheiderer's local-global principle (Proposition \ref{pro5.2} and Theorem \ref{thr5.3}), Scheiderer's Hessian criterion (Proposition \ref{pro5.6} and Theorem \ref{thr5.7}) and Marshall's boundary Hessian conditions (Proposition \ref{pro5.9} and Theorem \ref{thr5.10}). Similar to the matrix version of Schweighofer's Positivstellensatz given in section 4, we have a nice representation of diagonal polynomial matrices, but for an arbitrary symmetric polynomial matrices we need a 
denominator. 

\section{Preliminaries}

In this section we shall recall some basis concepts and  facts in Real algebraic geometry for matrices over  commutative rings which are  proposed by Schm\"{u}dgen (\cite{Schm-1}, \cite{Schm0}, \cite{Schm}) and Cimprič (\cite{Ci0}, \cite{Ci}).  

For $n\in \n^*$, let  $\hoa{M}_n(R)$  denote the ring of $n\times n$ matrices with entries from  a commutative unital ring  $R$. Denote by $\hoa{S}_n(R)$ the subset of  $\hoa{M}_n(R)$  consisting of all symmetric matrices. A subset $\hoa{M}$ of $\hoa{S}_n(R)$ is called a \textit{quadratic module}\footnotemark \footnotetext{In \cite{Schm-1} and \cite{Schm0}, the  term  \emph{m-admissible wedge} was used.} if 
$$ \i_n \in \hoa{M}, \quad \hoa{M} + \hoa{M} \subseteq \hoa{M}, \quad  A^T \hoa{M} A \subseteq \hoa{M}, \vm A\in \hoa{M}_n(R).$$
The smallest quadratic module which contains a given subset $\hoa{G}$ of $\hoa{S}_n(R)$ will be denoted by $\hoa{M}_{\hoa{G}}$. It is clear that 
$$ \hoa{M}_{\hoa{G}} = \{\sum_{i,j}A_{ij}^T G_i A_{ij}| G_i \in \hoa{G}\cup \{\i_n\}, A_{ij}\in \hoa{M}_n(R)\}. $$
In particular, a subset $M\subseteq R$ is a quadratic module if $1\in M, M+M\subseteq M, $ and $a^2M\subseteq M$ for all $a\in R$. The smallest quadratic module of $R$ which contains a given subset $G\subseteq  R$ will be denoted by $M_G$, and it consists of all  finite sums  of the form $\sum_{i,j} a_{ij}^2g_i$, $g_i\in G, a_{ij}\in R$.

A subset $\hoa{T}$ of $\hoa{S}_n(R)$ is called a \textit{preordering} if $\hoa{T}$ is a quadratic module in $\hoa{M}_n(R)$ and the set $\hoa{T}\cap (R\cdot \i_n)$ is closed under multiplication. The smallest preordering which contains a given subset $\hoa{G}$ of $\hoa{S}_n(R)$ will be denoted by $\hoa{T}_{\hoa{G}}$. We have 

\lm[{\cite[Lemma 2]{Ci}}] \label{lm1.1} For every subset $\hoa{G} $ of $\hoa{S}_n(R)$, 
$$\hoa{T}_{\hoa{G}} = \hoa{M}_{\hoa{G}\cup (\prod \hoa{G}'\cdot \i_n)},$$
where $\prod \hoa{G}' $ is the set of all finite product of elements from the set $\hoa{G}':=\{\bold{v}^T \bold{G} \bold{v}| \bold{G}\in \hoa{G}, \bold{v} \in R^n\}$. 
\elm 
In particular, a subset $T\subseteq R$ is a preordering  if  $T+T\subseteq T, T\cdot T \subseteq T, a^2\in T$ for every $a\in R$.   The smallest preordering of $R$ which contains a given subset $G\subseteq  R$ will be denoted by $T_G$. It is clear that 
$$ T_G=\{\sum_{\sigma=(\sigma_1,\cdots,\sigma_m) \in \{0,1\}^m}s_{\sigma} g_1^{\sigma_1}\cdots g_m^{\sigma_m}| m\in \n, g_i\in G, s_{\sigma} \in \sum R^2 \}, $$
where $\sum R^2$ is the set of all sums of squares of  finite elements from $R$. 

In the case $\hoa{G}=\emptyset$, $\sum_n R:=\hoa{M}_\emptyset = \hoa{T}_{\emptyset}$ is the set of all finite sums of elements of the form $A^TA$, where $A\in \hoa{M}_n(R)$, and which is the smallest quadratic module in $\hoa{M}_n(R)$.

For a quadratic module $M$ in $R$, denote 
$$ M^n:=\{\sum_{i} m_i A_i^T A_i | m_i\in M, A_i\in \hoa{M}_n(R)\}. $$
Then $M^n$ is the smallest quadratic module in $\hoa{M}_n(R)$ whose intersection with  $R\cdot \i_n$ is equal to $M\cdot \i_n$ (\cite[Proposition 3]{Ci}). 

\rem \label{rem1.2} \rm Let $M$ be a quadratic module of $R$. Denote by $D(d_1,\cdots,d_r)$, $r\leq n$, the $n\times n$ diagonal matrix 
 with  diagonal entries $d_1, \cdots, d_r, 0, \cdots, 0$, where  $d_i \in M $ for every $i=1,\cdots, r$. Then $D(d_1,\cdots, d_r) \in M^{n}$. \\
In fact, for every $i,j=1,\cdots, n$,  let $\bold{E}_{ij}$ be the coordinate matrices in $\hoa{M}_n(R)$. Note that for each $i=1,\cdots,n$, we have $ \bold{E}_{ii}= \bold{E}_{ii}^{T} \bold{E}_{ii}. $ Hence
$$ D(d_1,\cdots,d_r)=\sum_{i=1}^{r} d_i \bold{E}_{ii} = \sum_{i=1}^{r} d_i \bold{E}_{ii}^{T} \bold{E}_{ii} \in M^{n}. $$
\erem 

For any matrix $\bold{A}\in \hoa{M}_n(R)$, the  notation $\bold {A\succcurlyeq 0}$ means $\bold{A}$ is \emph{positive semidefinite},
 i.e. $\bold{x}^{T}\bold{A}\bold{x} \geq 0$ for every $\bold{x} \in \thuc^{n}$, and $\bold {A\succ 0}$ means $\bold{A}$ is \emph{positive definite}, i.e. $\bold{x}^{T}\bold{A}\bold{x} > 0$ for every $\bold{x} \in \thuc^{n}\tru \{0\}$.

In the following we consider $R$ to be the ring  $\thuc[X]:=\thuc[x_1,\cdots,x_d]$  of polynomials in $d$ variables $x_1,\cdots, x_d$ with real coefficients.  Then each element $\bold{A} \in  \hoa{M}_n(\thuc[X])$ is a matrix whose entries are polynomials from $\thuc[X]$, called a \emph{polynomial matrix}. Each element $\bold{A}\in \hoa{M}_n(\thuc[X])$ is also called a \emph{matrix polynomial},  because it can be viewed as a polynomial in $x_1,\cdots, x_d$ whose entries from $\hoa{M}_n(\thuc[X])$. Namely, we can write $\bold{A}$ as 
$$  \bold{A}=\sum_{|\alpha|=0}^{N} \bold{A}_{\alpha}X^{\alpha}, $$
where $\alpha=(\alpha_1,\cdots,\alpha_d) \in \n^{d}$, $|\alpha|:=\alpha_1+\cdots + \alpha_d$, $X^{\alpha}:=x_1^{\alpha_1}\cdots x_d^{\alpha_d}$, $\bold{A}_\alpha \in \hoa{M}_n(\thuc[X])$, $N$ is the maximum over all degree of entries of $\bold{A}$.  To unify notation, throughout the paper each element of $\hoa{M}_n(\thuc[X])$ is called a \emph{polynomial matrix}. 

   To every $\hoa{G}\subseteq \hoa{S}_n(\thuc[X])$ we associate the set 
$$ \bold{K}_{\hoa{G}}:=\{x\in \thuc^{d}| \bold{G}(x) \bold{\succcurlyeq 0}, \vm \bold{G}\in \hoa{G}\}. $$
In particular, for a subset $G$ of $\thuc[X]$, 
$$ K_G=\{x\in \thuc^d| g(x) \geq 0, \vm g\in G\}. $$
The following result of Cimprič (\cite{Ci}) shows that the set $\bold{K}_{\hoa{G}}$ can be determined by \emph{scalars}, i.e. by polynomials in $\thuc[X]$. 

\lm[{\cite[Proposition 5]{Ci}}] \label{lm1.3} Let $\hoa{G}\subset \hoa{S}_n(\thuc[X])$. Then there exists a subset $G$ of $ \thuc[X] $ with the following properties:
\it 
\item[(1)] $\bold{K}_{\hoa{G}} = K_{G}$;
\item[(2)] $(M_G)^{n} \subseteq \hoa{M}_{\hoa{G}}$;
\item[(3)] $(T_{G})^{n} \subseteq \hoa{T}_{\hoa{G}}$.
\hit
Moreover, if $\hoa{G}$ is finite then $G$ can be chosen to be finite. 
\elm 

It is well-known that every symmetric scalar matrix $A\in \hoa{S}_n(\thuc) $ can be diagonalized by an  orthogonal matrix $\bold{O}\in \hoa{M}_n(\thuc)$. For a polynomial matrix $\bold{A}$  in $\hoa{S}_n(\thuc[X])$, it is in general no longer true,  because   the matrix $\bold{O}$     may have rational entries (quotients of two polynomials in $\thuc[X]$).     However, Schm\"{u}dgen (\cite{Schm}) has showed that every symmetric polynomial matrix can be diagonalized by an invertible matrix in $\hoa{M}_n(\thuc[X])$ with a quotient by a non-zero polynomial in $\thuc[X]$. Moreover, in some special cases (e.g. that symmetric polynomial is in \textit{standard form}), that invertible matrix can be chosen to be lower triangular. 

\lm[{\cite[Corollary 9]{Schm}}] \label{lm1.4} Let $\bold{A}\in \hoa{S}_n(\thuc[X])$. Then there exist non-zero polynomials $b, d_j \in \thuc[X]$, $j=1,\cdots,r$, $r\leq n$, and matrices $\bold{X}_{+}, \bold{X}_{-} \in \hoa{M}_n(\thuc[X])$ such that 
$$ \bold{X}_{+}\bold{X}_{-} = \bold{X}_{-}\bold{X}_{+} = b\bold{I}_n, \quad b^2\bold{A} = \bold{X}_{+}\bold{D}\bold{X}_{+}^T, \quad \bold{D} = \bold{X}_{-} \bold{A} \bold{X}_{-}^T, $$
where $\bold{D}$ is the $n\times n$ diagonal matrix $D(d_1,\cdots,d_r)$.
\elm 
This lemma deduces a matrix version of  the theorem of Artin on Hilbert's seventeenth problem.

\coro[{\cite[Proposition 10]{Schm}}, \cite{GR}] \label{coro1.5} Let $\bold{F}\in \hoa{S}_n(\thuc[X])$. Then the following are equivalent:
\it 
\item[(1)] $\bold{F}\succcurlyeq \bold{0}$ (i.e. $\bold{F}(x) \succcurlyeq \bold{0}$ for every $x\in \thuc^d$);
\item[(2)] $b^2\bold{F} \in \sum_n \thuc[X]$ for some non-zero polynomial $b\in \thuc[X]$. 
\hit 
\ecoro

\section{Krivine-Stengle's Positivstellensatz for polynomial matrices}
In this section we shall  give a matrix version of Krivine-Stengle's  Positivstellensatz (cf. \cite{Kr}, \cite{St}, \cite[Positivstellensatz 2.2.1]{Ma}).  Let $\hoa{G}=\{\bold{G}_1,\cdots, \bold{G}_m\} \subseteq \hoa{S}_n(\thuc[X])$. Then by Lemma \ref{lm1.3}, there exists a subset $G=\{g_1,\cdots,g_k\}$ of $\thuc[X]$ such that $\bold{K}_{\hoa{G}} = K_G$ and $(T_G)^n \subseteq \hoa{T}_{\hoa{G}}$. 
For   diagonal polynomial matrices, we have the following 

\pro \label{pro2.1} Let $\bold{D}=D(d_1,\cdots,d_r)$, $r\leq n$, be an  $n\times n$ diagonal matrix in $\hoa{S}_n(\thuc[X])$. Then 
\it 
\item[(1)] $\bold{D} \succ \bold{0}$ on $\bold{K}_{\hoa{G}}$ if and only if there exist diagonal matrices $\bold{S}$ and $\bold{T}$ whose entries are in $T_G$ such that $\bold{S}\bold{D} = \bold{D} \bold{S} = \i_n + \bold{T}$.
\item[(2)] $\bold{D} \succcurlyeq  \bold{0}$ on $\bold{K}_{\hoa{G}}$ if and only if there exist an integer $m\geq 0$ and diagonal matrices $\bold{S}$ and $\bold{T}$ whose entries are in $T_G$ such that $\bold{S}\bold{D} = \bold{D} \bold{S} = \bold{D}^{2m} + \bold{T}$.
\item[(3)] $\bold{D} =  \bold{0}$ on $\bold{K}_{\hoa{G}}$ if and only if there exist an integer $m\geq 0$ such that $-\bold{D}^{2m} \in (T_G)^n$.
\item[(4)] $\bold{K}_{\hoa{G}} = \emptyset $ if and only if $-\i_n \in (T_G)^n$.
\hit 
\epro 
\pf  Note that in each of (1), (2), (3), (4), the "if" part is trivial. Therefore we shall prove the "only if" part in these statements. \\
(1) Assume $\bold{D} \succ \bold{0}$ on $\bold{K}_{\hoa{G}}$. Then  $r=n$ and $d_i > 0$ on $\bold{K}_{\hoa{G}} = K_G$ for all $i=1,\cdots, n$. It follows from Krivine-Stengle's  Positivstellensatz that for each $i=1,\cdots, n$, there exist $s_i$ and $t_i$ in $T_G$ such that $s_id_i=1 + t_i$. Then the matrices $\bold{S}=D(s_1,\cdots,s_r)  $ and $\bold{T}=D(t_1,\cdots,t_r)$  satisfy (1).\\
(2) Assume $\bold{D} \succcurlyeq  \bold{0}$ on $\bold{K}_{\hoa{G}}$. Then $d_i \geq  0$ on $\bold{K}_{\hoa{G}} = K_G$ for all $i=1,\cdots, r$. It follows from Krivine-Stengle's  Positivstellensatz that for each $i=1,\cdots, r$, there exist an  integer $m_i\geq 0$ and elements $s_i$ and $t_i$ in $T_G$ such that $s_id_i=d_i^{2m_i}+ t_i$.  Let $m=\max\{m_i, i=1,\cdots, r\}$.  Then  for every $i=1,\cdots, r$, we have 
$$ (s_id_i^{2(m-m_i)}) d_i=d_i^{2m}+ (t_id_i^{2(m-m_i)}). $$
Denote ${s'}_i:=s_id_i^{2(m-m_i)}$, ${t'}_i:=t_id_i^{2(m-m_i)}$. Then 
$\bold{S}=D({s'}_1,\cdots,{s'}_r)  $ and $\bold{T}=D({t'}_1,\cdots,{t'}_r)$  satisfy (2).\\
(3) Assume $\bold{D}= \bold{0}$ on $\bold{K}_{\hoa{G}}$. Then $d_i =  0$ on $\bold{K}_{\hoa{G}} = K_G$ for all $i=1,\cdots, r$. It follows from Krivine-Stengle's  Positivstellensatz that for each $i=1,\cdots, r$, there exists an  integer $m_i\geq 0$ such that $-d_i^{2m_i} \in T_G$. Then for $m=\max\{ m_i, i=1,\cdots, r\}$ we have $-d_i^{2m} \in T_G$ for every $i=1,\cdots, r$. Then $-D^{2m} \in (T_G)^n$ by Remark \ref{rem1.2}.\\
(4) follows from Krivine-Stengle' Positivstellensatz and Remark \ref{rem1.2}.
\epf 

For arbitrary symmetric polynomial matrices, we have the following 

\thr \label{thr2.2} Let  $\hoa{G}\subseteq \hoa{S}_n(\thuc[X])$, $G\subseteq \thuc[X]$, $\bold{K}_{\hoa{G}}$, $K_G$, $\hoa{T}_{\hoa{G}}$  and $T_G$ be determined as above. Then for $\bold{F} \in \hoa{S}_n(\thuc[X])$, we have 
\it 
\item[(1)] $\bold{F} \succ \bold{0}$ on $\bold{K}_{\hoa{G}}$ if and only if there exist a matrix $\bold{X}_{-}\in \hoa{M}_n(\thuc[X])$ and diagonal matrices $\bold{S}$ and $\bold{T}$ whose entries are in $T_G$ such that $\bold{S}(\bold{X}_{-}\bold{F}\bold{X}_{-}^T) = (\bold{X}_{-}\bold{F}\bold{X}_{-}^T) \bold{S} =   \i_n + \bold{T}$.
\item[(2)] $\bold{D} \succcurlyeq  \bold{0}$ on $\bold{K}_{\hoa{G}}$ if and only if there exist an integer $m\geq 0$, a matrix $\bold{X}_{-}\in \hoa{M}_n(\thuc[X])$ and diagonal matrices $\bold{S}$ and $\bold{T}$ whose entries are in $T_G$ such that $\bold{S}(\bold{X}_{-}\bold{F}\bold{X}_{-}^T) = (\bold{X}_{-}\bold{F}\bold{X}_{-}^T) \bold{S} = \bold{D}^{2m} + \bold{T}$.
\item[(3)] $\bold{D} =  \bold{0}$ on $\bold{K}_{\hoa{G}}$ if and only if there exist an integer $m\geq 0$  and  a matrix $\bold{X}_{-}\in \hoa{M}_n(\thuc[X])$ such that $-(\bold{X}_{-}\bold{F}\bold{X}_{-}^T)^{2m} \in (T_G)^n$.
\hit 
\ethr 
\pf  By Lemma \ref{lm1.4}, there exist non-zero polynomials $b, d_j \in \thuc[X]$, $j=1,\cdots,r$, $r\leq n$, and a matrix  $\bold{X}_{-} \in \hoa{M}_n(\thuc[X])$ such that  $ \bold{X}_{-} \bold{F} \bold{X}_{-}^T = D(d_1,\cdots, d_r)=:\bold{D}$. Note that 
$\bold{F} \succ \bold{0}$ (resp. $\succcurlyeq \bold{0}$, $=\bold{0}$) if and only if $\bold{D} \succ \bold{0}$ (and $r=n$) (resp. $\succcurlyeq \bold{0}$, $=\bold{0}$). Therefore the theorem follows from Proposition \ref{pro2.1}, applying for $\bold{D} = \bold{X}_{-} \bold{F} \bold{X}_{-}^T$.
\epf 

\rem \label{rem2.3} \rm 
\it 
\item[(1)] Theorem \ref{thr2.2} (1) gives a simpler representation of the positive definite polynomial matrix $\bold{F}$ on $\bold{K}_{\hoa{G}}$, comparing to the non-commutative version of Krivine-Stengle' Positivstellensatz given in \cite[sections  4.2 and 4.4]{Schm}. 
 \item[(2)] In \cite{Ci} the author has given a matrix version of Krivine-Stengle's Positivstellensatz without the matrix $\bold{X}_{-}$ in representation of $\bold{F}$.  He requires also   $\bold{S}$, $\bold{T} \in (T_G)^n$, however they are in general not diagonal. 
\hit 
\erem

\section{Schweighofer's Positivstellensatz  for polynomial matrices}
In this section we give a matrix version of Schweighofer's Positivstellensatz (\cite{Schw}) which is recalled as follows. For a polynomial $f\in \thuc[X]$ and a subset $S\subseteq \thuc^d$,   a real number $y\in \thuc$ is called an \textit{asymptotic value } of $f$ on $S$ if there exists a sequence $(x_k)_{k\in \n}\subseteq S$ such that $\dlim_{k\mtn \vc} ||x_k|| = \vc$  and $\dlim_{k\mtn \vc } f(x_k) = y$. Denote by $R_{\vc}(f,S)$ the set of all asymptotic values of $f$ on $S$. Then we have 

\thr[{\cite[Theorem 9]{Schw}}] \label{thr3.1}  Let $G=\{g_1,\cdots, g_m\} \subseteq \thuc[X]$, and $f\in \thuc[X]$. Assume
\it 
\item[(1)] $f>0$ on $K_G$;
\item[(2)] $f$ is bounded on $K_G$;
\item[(3)] $R_{\vc}(f, K_G)$ is a finite subset of $ \thuc_{+}$.
\hit
Then $f\in T_G$.
\ethr 

We give firstly a version of this theorem for diagonal polynomial matrices. 

\pro \label{pro3.2} Let $\hoa{G}=\{\bold{G}_1,\cdots, \bold{G}_m\} \subseteq \hoa{S}_n(\thuc[X])$ and $\bold{D}=D(d_1,\cdots,d_n)$ be an $n \times n$ diagonal matrix  in $\hoa{S}_n(\thuc[X])$ with $d_i \not = 0$ for every $i=1,\cdots, n$.  Assume 
\it 
\item[(1)] $\bold{D}\succ \bold{0}$ on $\bold{K}_{\hoa{G}}$;
\item[(2)] $\bold{D}$ is bounded on $\bold{K}_{\hoa{G}}$ (i.e. there exists a number  $N\in \thuc_+$ such that $N \cdot  \i_n \pm D \succcurlyeq \bold{0}$ on $\bold{K}_{\hoa{G}}$); 
\item[(3)] For every $i=1,\cdots, n$, $R_\vc(d_i, \bold{K}_{\hoa{G}}) $ is a finite subset of $ \thuc_{+}$.
\hit 
Then there exists a finite subset $G$ of $~\thuc[X]$ such that $D\in (T_G)^n \subseteq \hoa{T}_{\hoa{G}}$. 
 \epro 

\pf  By Lemma \ref{lm1.3}, there exists a finite subset $G$ of $ \thuc[X]$ such that $\hoa{K}_{\hoa{G}} = K_G$ and $(T_G)^n \subseteq \hoa{T}_{\hoa{G}}$. By hypothesis, for every $i=1,\cdots, n$ we have 
\it 
\item $d_i > 0$ on $K_G$; 
\item $d_i$ is bounded on $K_G$;
\item $R_\vc(d_i, K_G)$ is a finite subset of $\thuc_{+}$.
\hit
Then it follows from Theorem \ref{thr3.1} that $d_i\in T_G$ for every $i=1,\cdots, n$. This implies that $\bold{D} \in (T_G)^n$ by Remark \ref{rem1.2}.
\epf 
For arbitrary symmetric polynomial matrices we have the following 
\thr \label{thr3.3} Let $\hoa{G}=\{\bold{G}_1,\cdots, \bold{G}_m\} \subseteq \hoa{S}_n(\thuc[X])$ and $\bold{F} \in \hoa{S}_n(\thuc[X])$. Assume 
\it 
\item[(1)] $\bold{F} \succ \bold{0}$ on $\hoa{K}_{\hoa{G}}$;
\item[(2)]  $\bold{F}$ is bounded on $\bold{K}_{\hoa{G}}$ (i.e. there exists a number $N\in \thuc_+ $ such that $N\cdot \i_n \pm \bold{F} \succcurlyeq \bold{0}$ on $\bold{K}_{\hoa{G}}$);
\item[(3)] for every $\bold{x}\in \thuc^{n}\tru \{0\}$, $R_\vc(\bold{x}^{T}\bold{F}\bold{x}, \bold{K}_{\hoa{G}})$ is a finite subset of $\thuc_+$.
\hit 
Then there exist a finite subset $G$ of $\thuc[X]$ and 
\it 
\item[(i)] a matrix $\bold{X}_{-} \in \hoa{M}_n(\thuc[X])$ such that $\bold{X}_{-}\bold{F}\bold{X}_{-}^T \in(T_G)^n \subseteq \hoa{T}_{\hoa{G}}$;
\item[(ii)] a non-zero polynomial $b\in \thuc[X]$   such that $b^2\bold{F} \in (T_G)^n \subseteq \hoa{T}_{\hoa{G}}$.
\hit
\ethr 
\pf  By Lemma \ref{lm1.4},  there exist non-zero polynomials $b, d_j \in \thuc[X]$, $j=1,\cdots,r$, $r\leq n$, and matrices $\bold{X}_{+}, \bold{X}_{-} \in \hoa{M}_n(\thuc[X])$ such that 
$$ \bold{X}_{+}\bold{X}_{-} = \bold{X}_{-}\bold{X}_{+} = b\bold{I}_n, \quad b^2\bold{F} = \bold{X}_{+}\bold{D}\bold{X}_{+}^T, \quad \bold{D} = \bold{X}_{-} \bold{F} \bold{X}_{-}^T,$$
where $\bold{D}=D(d_1,\cdots,d_r)$ is the $n\times n$ diagonal polynomial matrix. Since $\bold{F}\succ \bold{0}$ on $\bold{K}_{\hoa{G}}$, $r=n$. Note that for every $i=1,\cdots, n$, 
\begin{equation} \label{equ3.1.1}
 d_i = \bold{e}_i^{T} \bold{D} \bold{e}_i = ({\bold{X}_{-}}^{T}\bold{e}_i)^{T} \bold{F} ({\bold{X}_{-}}^{T}\bold{e}_i), 
\end{equation} 
where $\bold{e}_i$, $i=1,\cdots, n$, are the coordinate vectors in $\thuc^{n}$. Since $\bold{v}_i:={\bold{X}_{-}}^{T}\bold{e}_i \in \thuc^{n}\tru \{0\}$ and $\bold{F}\succ \bold{0}$ on $\bold{K}_{\hoa{G}}$, it follows that $d_i>0$ on $\bold{K}_{\hoa{G}}$ for every $i=1,\cdots,n$.

By (2) and in view of (\ref{equ3.1.1}), for each $i=1,\cdots, n$, we have 
$$ N (\bold{v}_i^{T}\bold{v}_i) \pm d_i \geq 0 \mbox{ on } \bold{K}_{\hoa{G}}.$$
It follows that each $d_i$, $i=1,\cdots, n$, is bounded on $\bold{K}_{\hoa{G}}$. Moreover, by (3) and in view of (\ref{equ3.1.1}), $R_\vc(d_i, \bold{K}_{\hoa{G}})$ is a finite subset of $\thuc_+$ for each $i=1,\cdots, n$. Then it follows from Proposition \ref{pro3.2} that there exists a finite subset $G$ of $\thuc[X]$ such that $\bold{D} \in (T_G)^{n} \subseteq \hoa{T}_{\hoa{G}}$, hence $\bold{X}_{-}\bold{F} \bold{X}_{-}^T \in (T_G)^n$, i.e. we have (i). Moreover, since $(T_G)^{n}$ is a quadratic module of $\hoa{M}_n(\thuc[X])$, by definition we have $b^{2}\bold{F} = \bold{X}_{+}\bold{D} \bold{X}_{+}^{T} \in (T_G)^{n}$, i.e. we have (ii). The proof is complete.
\epf 

\section{Archimedean quadratic modules}
In this section we deal with    Archimedean  quadratic modules of polynomial matrices.  We recall the definition  of Archimedean quadratic modules, and show  that the Archimedean property of a quadratic module $M$ in $\thuc[X]$ is the same as that of the quadratic module $M^n$ in $\hoa{M}_n(\thuc[X])$. Moreover, we also show that the compactness of $\bold{K}_{\hoa{G}}$ is equivalent to the Archimedean property of the preordering $\hoa{T}_{\hoa{G}}$.

\df[\cite{Schm0}, \cite{Schm}, \cite{Ci0}] \label{df4.1} \rm Let $\hoa{M}$ be a quadratic module  in $ \hoa{M}_n(\thuc[X])$. 
\it  
\item[(1)]  $\hoa{M}$  is called \emph{Archimedean} if for each element $\bold{A} \in \hoa{M}_n(\thuc[X])$ there exists a number $n \in \n$ such that $n\cdot \i_n - \bold{A}^{T}\bold{A} \in \hoa{M}.$
\item[(2)] Denote 
$$ H_{\hoa{M}}:=H_{\hoa{M}}(\hoa{M}_n(\thuc[X])):= \{\hoa{A}\in \hoa{M}_n(\thuc[X]) | \tt r\in \thuc_+: r^2\cdot \i_n - \hoa{A}^{T}\hoa{A} \in \hoa{M}\}. $$
\hit 
It is  clear  that  the quadratic module $\hoa{M}$ in $ \hoa{M}_n(\thuc[X])$ is Archimedean if and only if  $H_{\hoa{M}} = \hoa{M}_n(\thuc[X])$. \footnotemark \footnotetext{In this case,  the ring $\hoa{M}_n(\thuc[X])$ is called \emph{algebraically bounded} with respect to the quadratic module $\hoa{M}$, cf. \cite{Schm0}.}  Moreover, $H_{\hoa{M}}$ is a subring of $\hoa{M}_n(\thuc[X])$ (cf. \cite[Corollary 2.2]{Schm0}, \cite[Corollary 5]{Ci0}), and  it is called the  \textit{  ring of bounded elements } of $\hoa{M}_n(\thuc[X])$ with respect to the quadratic module $\hoa{M}$.
\edf 
The following fact is useful and it is easy to check (cf. \cite[Lemma 2.1(ii)]{Schm0}, \cite[Lemma 3]{Ci0}).
\lm  \label{lm4.3}  Let $\hoa{M}$ be a quadratic module  in $ \hoa{M}_n(\thuc[X])$. Then for any $\bold{A}\in \hoa{S}_n(\thuc[X])$ and for any $r\in \thuc_+$, we have $r^2 \cdot \i_n- \bold{A}^2 \in \hoa{M}$ if and only if $r \cdot \i_n \pm \bold{A} \in \hoa{M}$.
\elm 

Similar to the case of polynomials (cf. \cite[Corollary 5.2.4]{Ma}), we can check the Archimedean property of quadratic modules  of polynomial matrices simply  as follows.

\pro \label{pro4.4} Let $\hoa{M}\subseteq \hoa{M}_{n}(\thuc[X])$ be a quadratic module. Then the following are equivalent:
\it 
\item[(1)] $\hoa{M}$ is Archimedean.
\item[(2)] $r\cdot \i_n -\sum_{i=1}^d x_i^2\cdot \i_n \in \hoa{M}$ for some positive real number $r$.
\item[(3)] $r\cdot \i_n \pm x_i\cdot \i_n \in \hoa{M}$ for some positive real number $r$.
\hit
\epro 
\pf 
$(1) \Sr (2)$ is clear. If (2) holds, for each $i=1,\cdots, d$ we have 
$$ r\cdot \i_n - x_i^2\cdot \i_n = (r \cdot \i_n - \sum_{i=1}^dx_i^2 \cdot \i_n) + \sum_{j\not = i} x_j^2\cdot \i_n \in \hoa{M}.  $$ 
It follows from Lemma  \ref{lm4.3}  that $\can{r}\cdot \i_n\pm x_i\cdot \i_n \in \hoa{M}$ for every $i=1,\cdots, d$, i.e. we have (3). \\
To show  $(3) \Sr (1)$, it suffices to prove  that $H_{\hoa{M}} = \hoa{M}_{n}(\thuc[X])$.  Since $\hoa{M}_n(\thuc[X])$ is generated as an $\thuc$-algebra by $x_i, i=1,\cdots, d$, and the  coordinate matrices $\bold{E}_{ij}, i, j =1,\cdots, n$, of $M_n(\thuc)$, and since $H_{\hoa{M}}$ is closed under addition and  multiplication, it is enough to show that $x_i\cdot \i_n \in H_{\hoa{M}}$ for every $i=1,\cdots, d$ and $\bold{E}_{ij}\in H_{\hoa{M}}$ for every $i,j=1,\cdots,n$. \\
Since  $r\cdot \i_n \pm x_i\cdot \i_n \in \hoa{M}$, it follows from  Lemma \ref{lm4.3} that $r^2\cdot \i_n - x_i^2\cdot \i_n \in \hoa{M}$, hence $x_i\cdot \i_n \in H_{\hoa{M}}$ for every $i=1,\cdots, d$. On the other hand, for every $i,j=1,\cdots, n$, we have 
$\bold{E}_{ij}^{T} \bold{E}_{ij} = \bold{E}_{jj}$. Therefore, 
$$ \i_n - \bold{E}_{ij}^T \bold{E}_{ij} = \sum_{k\not = j}\bold{E}_{kk}^T \bold{E}_{kk} \in \hoa{M}. $$
It follows that $\bold{E}_{ij} \in H_{\hoa{M}}$ for every $i,j=1,\cdots, n$. The proof is complete.
\epf 
Using this criterion we can show now the equivalence  of the Archimedean property of a quadratic module $M$ in $\thuc[X]$ and the quadratic module $M^n$ in $\hoa{M}_n(\thuc[X])$.
\pro \label{pro4.5} Let $M$ be a quadratic module in $\thuc[X]$. Then $M$ is Archimedean if and only if $M^n$ is an Archimedean quadratic module in $\hoa{M}_n(\thuc[X])$.
\epro 
\pf  The "only if"  part follows easily from the usual criterion for Archimedean property of quadratic modules in $\thuc[X]$   (cf. \cite[Corollary 5.2.4]{Ma}) and Proposition \ref{pro4.4}. Now we prove  the "if" part.\\
Assume $M^n$ is Archimedean. Then it follows from Proposition \ref{pro4.4} that  for every $i=1,\cdots, d$, we have   $r\cdot \i_n \pm x_i\cdot \i_n\in M^n  $ for some $r\in \thuc_{+}$. Then we can write 
$$ (r\pm x_i)\cdot \i_n= r\cdot \i_n \pm x_i\cdot \i_n = \sum_{j=1}^m m_j \bold{A}_j^T \bold{A}_j, \mbox{where } m_j\in M, \bold{A}_j\in \hoa{M}_n(\thuc[X]).$$
Note that $Tr(\bold{A}_j^T \bold{A}_j) \in \sum \thuc[X]^2$ for each $j=1,\cdots, m$. Then  for every $i=1,\cdots, d$ we have
$$ r \pm x_i   = \ps{1}{n}Tr\big((r\pm x_i)\cdot \i_n\big) = \ps{1}{n} \sum_{j=1}^m m_j Tr(\bold{A}_j^T\bold{A}_j) \in M.$$
Hence $M$ is Archimedean (cf. \cite[Corollary 5.2.4]{Ma}).
\epf 
It is well-known that the compactness of the basic semi-algebraic set $K_G\subseteq \thuc^n$, $G\subseteq \thuc[X]$, is equivalent to the Archimedean property of the preordering $T_G$ in $\thuc[X]$ (cf. \cite[Theorem 6.1.1]{Ma}).  For polynomial matrices we have also the same result.

\pro \label{pro4.6}  Let $\hoa{G} \subseteq \hoa{S}_n(\thuc[X])$. Then $\bold{K}_{\hoa{G}}$ is compact if and only if $\hoa{T}_{\hoa{G}}$ is Archimedean.
\epro 
\pf  Assume $\hoa{T}_{\hoa{G}}$ is Archimedean. It follows from Proposition \ref{pro4.4} that there exists a number $r\in \thuc_+$ such that $r\cdot \i_n -\sum_{i=1}^d x_i^2 \cdot \i_n  \in \hoa{T}_{\hoa{G}}$. This implies $r\cdot \i_n -\sum_{i=1}^d x_i^2 \cdot \i_n  \succcurlyeq \bold{0}$ on $\bold{K}_{\hoa{G}}$. Then for any point $p=(p_1,\cdots,p_d) \in \bold{K}_{\hoa{G}}$, we have 
$r- \sum_{i=1}^d p_i^2 \geq 0$, i.e., $|| p|| \leq \can{r}$. It follows that $\bold{K}_{\hoa{G}}$ is bounded, whence compact. 

Conversely, assume that $\bold{K}_{\hoa{G}}$ is compact. By Lemma \ref{lm1.3}, there exists a subset $G$ of $\thuc[X]$ such that $\bold{K}_{\hoa{G}} = K_G$ and $(T_G)^n \subseteq \hoa{T}_{\hoa{G}}$. Then $K_G$ is compact. It follows that $T_G$ is  an Archimedean  quadratic module in $\thuc[X]$ (cf. \cite[Theorem 6.1.1]{Ma}). Then $(T_G)^n \subseteq \hoa{M}_n(\thuc[X])$ is Archimedean by Proposition \ref{pro4.5}. This implies that $\hoa{T}_{\hoa{G}} \supseteq (T_G)^n$ is Archimedean.
\epf 

\rem \label{rem4.7}\rm  For any  $\hoa{G}\subseteq \hoa{S}_n(\thuc[X])$,  since $\hoa{M}_{\hoa{G}} \subseteq \hoa{T}_{\hoa{G}}$, if $\hoa{M}_{\hoa{G}}$ is Archimedean then $ \hoa{T}_{\hoa{G}}$ is Archimedean, hence $\bold{K}_{\hoa{G}}$ is compact by 
Proposition \ref{pro4.6}. The converse is in general not true, even for polynomials (i.e. for $n=1$). A natural question, like Putinar's  question for polynomials (cf. \cite[Chapter 7]{Ma}), is that in which cases the compactness of $\bold{K}_{\hoa{G}}$ implies the 
Archimedean property of $\hoa{M}_{\hoa{G}}$?  For univariate polynomial matrices, we have a  confirmation.
\erem 

\pro \label{pro4.8} Let $\thuc[t]$ be the ring of polynomial in one variable $t$ with real coefficients. Then, for a finite set  $\hoa{G}\subseteq \hoa{S}_n(\thuc[t])$,  if  $\bold{K}_{\hoa{G}}$ is compact then $\hoa{M}_{\hoa{G}}$ is Archimedean.
\epro 
\pf  By the same argument  as given in the proof of the "only if" part  of   Proposition \ref{pro4.6}, using  \cite[Theorem 7.1.2]{Ma} instead of  \cite[Theorem 6.1.1]{Ma}, we obtain the result. 
\epf 
For  multivariate polynomial matrices (i.e. for $d\geq 2$), the compactness of $\bold{K}_{\hoa{G}}$ is in general not sufficient to deduce the Archimedean property of $\hoa{M}_{\hoa{G}}$. It is even not true for the case of multivariate polynomials (i.e. for $d\geq 2$ and $n=1$), see, for example,  Jacobi-Prestel's  counterexample (cf. \cite[Example 4.6]{JP}). 

\section{Local-global principle and Hessian conditions for polynomial matrices}

For a set $\hoa{G} $  and a polynomial matrix $\bold{F}$ in  $\hoa{S}_n(\thuc[X])$, it is obvious that if $\bold{F} \in \hoa{T}_{\hoa{G}}$ (resp. $\bold{F} \in \hoa{M}_{\hoa{G}}$) then $\bold{F} \succcurlyeq \bold{0}$ on $\bold{K}_{\hoa{G}}$.  The converse is true only in some special cases. For example, if $\bold{K}_{\hoa{G}}$ is compact (equivalently, $\hoa{T}_{\hoa{G}}$ is Archimedean by Proposition \ref{pro4.6}) (resp.  if $\hoa{M}_{\hoa{G}}$ is Archimedean) and $\bold{F}\succ \bold{0}$ on $\bold{K}_{\hoa{G}}$ then $\bold{F}\in \hoa{T}_{\hoa{G}}$ (resp. $\bold{F} \in \hoa{M}_{\hoa{G}}$). This is a matrix version of Schm\"{u}dgen's  
Positivstellensatz, see, for example \cite{CZ}  (resp. Putinar's  Positivstellensatz, see, for example  \cite{SchH} or \cite{KSch}). \\
 In  the case where $\bold{K}_{\hoa{G}}$ is not compact, we have given in section 4 some special  conditions for $\bold{F}\succ \bold{0}$  to ensure that  $\bold{X}_{-}\bold{F}\bold{X}_{-}^T$ or $b^{2}\bold{F}$ belongs to $\hoa{T}_{\hoa{G}}$.  If $\bold{F}$ vanishes at some points in $\bold{K}_{\hoa{G}}$, we need some conditions at these zeros to ensure for $\bold{F}$ belonging to $\hoa{T}_{\hoa{G}}$ or $\hoa{M}_{\hoa{G}}$. In  the  polynomial case (i.e. $n=1$), one of the well-known criterion for $f\geq 0$ on $K_G$ to be in  $T_G$ (resp. $M_G$) is the \textit{Hessian criterion of Scheiderer} (cf. \cite{Sch1}, \cite{Sch2} or \cite[section 9.5]{Ma}), and to prove it, he used his \textit{local-global principle } (cf.  \cite{Sch1} or \cite[section 9.2]{Ma}). Moreover, Marshall (\cite{Ma}) has given  \textit{boundary  Hessian conditions} for $f$  to ensure that  $f\in M_G$ whenever it is non-negative on $K_G$. Therefore,  in this section we give a matrix version of the local-global principle of Scheiderer, Scheiderer's Hessian criterion and the boundary 
Hessian conditions of Marshall.  
\subsection{Local-global principle for polynomial matrices}

First we recall Scheiderer's  local-global principle. 

\thr[{\cite{Sch1}, \cite[Theorem 9.2.1]{Ma}}] \label{thr5.1} Let $G=\{g_1,\cdots,g_m\} \subseteq \thuc[X]$ and $f\in \thuc[X]$. Assume 
\it  
\item[(1)] $K_G$ is compact;
\item[(2)] $f\geq 0$ on $K_G$, and $f$ has only finitely many zeros in $K_G$;
\item[(3)] at each zero $p$ of $f$ in $K_G$, $f\in (\widehat{T_G})_p \subseteq \thuc[[X-p]]$, the preordering of $\thuc[[X-p]]$ generated by $G$.
\hit
Then $f\in T_G$.
\ethr
For any subset $\hoa{G}=\{G_1,\cdots,G_m\}$ of $\hoa{S}_n(\thuc[X])$, by Lemma \ref{lm1.3}, there exists a finite subset $G$ of $\thuc[X]$ such that $\bold{K}_{\hoa{G}} = K_G$ and $(T_G)^n \subseteq \hoa{T}_{\hoa{G}}$. We firstly give a local-global principle for diagonal polynomial matrices.

\pro \label{pro5.2} Let $\hoa{G}=\{G_1,\cdots,G_m\} \subseteq \hoa{S}_n(\thuc[X])$ and $G\subseteq \thuc[X]$ as above. Let $\bold{D} = D(d_1,\cdots,d_r)$, $r\leq n$,  be an $n\times n$ diagonal polynomial matrix in $\hoa{S}_n(\thuc[X])$. Assume 
\it  
\item[(1)] $\bold{K}_{\hoa{G}}$ is compact;
\item[(2)] $\bold{D}\geq 0$ on $\bold{K}_{\hoa{G}}$, and each $d_i$ has only finitely many zeros in $\bold{K}_{\hoa{G}}$;
\item[(3)] at each zero $p$ of each $d_i$ in $\bold{K}_{\hoa{G}}$, $d_i\in (\widehat{T_G})_p \subseteq \thuc[[X-p]]$.
\hit
Then $\bold{D} \in (T_G)^n \subseteq \hoa{T}_{\hoa{G}}$.
\epro 
\pf Theorem \ref{thr5.1}, applying for each $d_i$, implies that each  $d_i$ belongs to $T_G$. Then $\bold{D} \in (T_G)^n \subseteq \hoa{T}_{\hoa{G}}$ by Remark \ref{rem1.2}.
\epf 
For arbitrary polynomial matrices, we have the following

\thr \label{thr5.3} Let $\hoa{G}=\{G_1,\cdots,G_m\} \subseteq \hoa{S}_n(\thuc[X])$ and $G\subseteq \thuc[X]$ as above. Let $\bold{F}\in \hoa{S}_n(\thuc[X])$. Assume
\it
\item[(1)] $\bold{K}_{\hoa{G}}$ is compact;
\item[(2)] $\bold{F}\succcurlyeq \bold{0}$ on $\bold{K}_{\hoa{G}}$;
\item[(3)] for each $\bold{x} \in \thuc^n\tru \{0\}$, $\bold{x}^T\bold{F}\bold{x} $ has only finitely many zeros in $\bold{K}_{\hoa{G}}$;
\item[(4)] for each $\bold{x} \in \thuc^n\tru \{0\}$ and for  each zero $p$ of $\bold{x}^T\bold{F}\bold{x} $ in $\bold{K}_{\hoa{G}}$, $\bold{x}^T\bold{F}\bold{x}$ belongs to $(\widehat{T_G})_p \subseteq \thuc[[X-p]]$.
\hit 
Then 
\it 
\item[(i)] there exists a matrix $\bold{X}_{-}\in \hoa{M}_{n}(\thuc[X])$ such that $\bold{X}_{-}\bold{F}\bold{X}_{-}^T \in (T_G)^n\subseteq \hoa{T}_{\hoa{G}}$;
\item[(ii)] there exists a non-zero polynomial $b\in \thuc[X]$ such that $b^2\bold{F} \in (T_G)^n \subseteq \hoa{T}_{\hoa{G}}$.
\hit
\ethr 
\pf 
By Lemma \ref{lm1.4},  there exist non-zero polynomials $b, d_j \in \thuc[X]$, $j=1,\cdots,r$, $r\leq n$, and matrices $\bold{X}_{+}, \bold{X}_{-} \in \hoa{M}_n(\thuc[X])$ such that 
$$ \bold{X}_{+}\bold{X}_{-} = \bold{X}_{-}\bold{X}_{+} = b\bold{I}_n, \quad b^2\bold{F} = \bold{X}_{+}\bold{D}\bold{X}_{+}^T, \quad \bold{D} = \bold{X}_{-} \bold{F} \bold{X}_{-}^T,$$
where $\bold{D}=D(d_1,\cdots,d_r)$ is the $n\times n$ diagonal polynomial matrix.  Note that for every $i=1,\cdots, r$, 
\begin{equation} \label{equ3.1}
 d_i = \bold{e}_i^{T} \bold{D} \bold{e}_i = ({\bold{X}_{-}}^{T}\bold{e}_i)^{T} \bold{F} ({\bold{X}_{-}}^{T}\bold{e}_i).
\end{equation} 
Since $\bold{v}_i:={\bold{X}_{-}}^{T}\bold{e}_i \in \thuc^{n}\tru \{0\}$ and $\bold{F}\succcurlyeq \bold{0}$ on $\bold{K}_{\hoa{G}}$, it follows that $d_i\geq 0$ on $\bold{K}_{\hoa{G}}$ for every $i=1,\cdots,r$.

By (3) and in view of (\ref{equ3.1}), each $d_i$ has only finitely many zeros in $K_G$. By (4) and in view of (\ref{equ3.1}), at each zero $p$ of $d_i$ in $K_G$, $d_i\in (\widehat{T_G})_p$. It follows from Proposition \ref{pro5.2} that 
$\bold{D} \in (T_G)^n$, hence $\bold{X}_{-}\bold{F} \bold{X}_{-}^T \in (T_G)^n$, i.e. we have (i). Moreover, since $(T_G)^{n}$ is a quadratic module of $\hoa{M}_n(\thuc[X])$, by definition we have $b^{2}\bold{F} = \bold{X}_{+}\bold{D} \bold{X}_{+}^{T} \in (T_G)^{n}$, i.e. we have (ii). The proof is complete.
\epf

\subsection{Hessian criterion for polynomial matrices} 
We recall  firstly Scheiderer's Hessian criterion for polynomials in $\thuc[X]$. 
\thr[{\cite[Example 3.18]{Sch1},\cite[Corollary 3.6]{Sch2}}] \label{thr5.4}  Let $G=\{g_1,\cdots,g_m\}$ be a subset of $ \thuc[X]$  and $f\in \thuc[X]$. Assume 
\it  
\item[(1)] $K_G$ is compact (resp. the quadratic module $M_G$  is Archimedean);
\item[(2)] $f\geq 0$ on $K_G$; 
\item[(3)] $f$ has only finitely many zeros in $K_G$ and all of them are in the interior of $K_G$;
\item[(4)] at each zero $p$ of $f$ in $K_G$, the Hessian  $D^2f(p)$ of $f$ at $p$  is positive definite.
\hit
Then $f\in T_G$ (resp. $f\in M_G$). 
\ethr
\rem \label{rem5.5} \rm
\it 
\item[(1)] Condition (3) in Theorem \ref{thr5.4} requires each zero $p$ of $f$ in $K_G$ must be in the interior of $K_G$, then it follows that $p$ is a  local minimum  of $f$ in $K_G$. Therefore, for  a Taylor expansion of $f$ in a neighborhood of $p$, $f=f_0+f_1+f_2+\cdots \in \thuc[[X-p]]$, we have $f_0=f_1=0$. Moreover, this condition implies that $(\widehat{T_G})_p = \sum \thuc[[X-p]]^2$. 
\item[(2)] The Hessian condition of $f$ at $p$ in Theorem \ref{thr5.4} implies that in the Taylor expansion  of $f$ in a neighborhood of $p$, the quadratic form $f_2$ can be written as $x_1^2 + \cdots + x_d^2$ (after changing coordinates). Then by some special techniques and using Local-global principle (Theorem \ref{thr5.1}), we have the conclusion for $T_G$.
 \hit
\erem
Like in previous sections, we fist give a result for diagonal polynomial matrices.

\pro \label{pro5.6} Let $\hoa{G}=\{G_1,\cdots,G_m\} \subseteq \hoa{S}_n(\thuc[X])$ and $G\subseteq \thuc[X]$ as in Lemma \ref{lm1.3}. 
Let $\bold{D} = D(d_1,\cdots,d_r)$, $r\leq n$,  be an $n\times n$ diagonal polynomial matrix in $\hoa{S}_n(\thuc[X])$. Assume 
\it  
\item[(1)] $\bold{K}_{\hoa{G}}$ is compact (resp. $M_G$ is Archimedean);
\item[(2)] $\bold{D}\geq 0$ on $\bold{K}_{\hoa{G}}$;
\item[(3)] each $d_i$ has only finitely many zeros in $\bold{K}_{\hoa{G}}$, and all of  them lie  in the interior of $\bold{K}_{\hoa{G}}$;
\item[(4)] at each zero $p$ of each $d_i$ in $\bold{K}_{\hoa{G}}$, the Hessian $D^2d_i(p)$ is positive definite.
\hit
Then $\bold{D} \in (T_G)^n \subseteq \hoa{T}_{\hoa{G}}$ (resp. $\bold{D} \in (M_G)^n \subseteq \hoa{M}_{\hoa{G}}$).
\epro
\pf Theorem \ref{thr5.4}, applying for each $d_i$, implies that each  $d_i$ belongs to $T_G$ (resp. $M_G$). Then 
$\bold{D} \in (T_G)^n$ (resp. $\bold{D} \in (M_G)^n$) by Remark \ref{rem1.2}. 
\epf 

For arbitrary polynomial matrices we have the following
\thr \label{thr5.7} Let $\hoa{G}=\{G_1,\cdots,G_m\} \subseteq \hoa{S}_n(\thuc[X])$ and $G\subseteq \thuc[X]$ as in Lemma \ref{lm1.3}. 
Let $\bold{F}\in \hoa{S}_n(\thuc[X])$. Assume
\it
\item[(1)] $\bold{K}_{\hoa{G}}$ is compact (resp. $M_G$ is Archimedean);
\item[(2)] $\bold{F}\succcurlyeq \bold{0}$ on $\bold{K}_{\hoa{G}}$;
\item[(3)] for each $\bold{x} \in \thuc^n\tru \{0\}$, $\bold{x}^T\bold{F}\bold{x} $ has only finitely many zeros in $\bold{K}_{\hoa{G}}$ and each zero lies in the interior of $\bold{K}_{\hoa{G}}$;
\item[(4)] for each $\bold{x} \in \thuc^n\tru \{0\}$ and for  each zero $p$ of $\bold{x}^T\bold{F}\bold{x} $ in $\bold{K}_{\hoa{G}}$, the Hessian $D^2 \big(\bold{x}^T\bold{F}\bold{x}\big)(p)$ is positive definite.
\hit 
Then 
\it 
\item[(i)] there exists a matrix $\bold{X}_{-}\in \hoa{M}_{n}(\thuc[X])$ such that $\bold{X}_{-}\bold{F}\bold{X}_{-}^T \in (T_G)^n\subseteq \hoa{T}_{\hoa{G}}$ (resp.  $\bold{X}_{-}\bold{F}\bold{X}_{-}^T \in (M_G)^n\subseteq \hoa{M}_{\hoa{G}}$);
\item[(ii)] there exists a non-zero polynomial $b\in \thuc[X]$ such that $b^2\bold{F} \in (T_G)^n \subseteq \hoa{T}_{\hoa{G}}$ (resp. $b^2\bold{F} \in (M_G)^n \subseteq \hoa{M}_{\hoa{G}}$).
\hit
\ethr 
\pf  By a similar argument to the one given  in the proof of  Theorem \ref{thr5.3}, using Proposition \ref{pro5.6}, we have the proof.
\epf

\subsection{Boundary Hessian conditions for polynomial matrices}
Let us recall the  boundary Hessian conditions of a polynomial at a point, which is defined by Marshall (cf. \cite[section 9.5]{Ma}). Let $G \subseteq \thuc[X]$ and $f\in \thuc[X]$. We say that $f$ satisfies the \textit{boundary Hessian conditions (BHC)} at a point $p\in K_G$ with respect to $t_1,\cdots, t_k$, $1\leq k\leq d$, which are part of a system of uniformizing parameters $t_1,\cdots, t_d$ at $p$,   if $p$ is a non-singular point of $\thuc^d$, and in the completion $\thuc[[t_1,\cdots,t_d]]$ of $\thuc[X]$ at $p$, $f$ decomposes as $f=f_0+f_1+f_2+\cdots$ (where $f_j$ is homogeneous of degree $j$ in  the variables $t_1,\cdots,t_d$ with coefficients in $\thuc$), $f_1=a_1t_1+\cdots + a_kt_k$, $a_i>0$ for $i=1,\cdots, k$, and the quadratic form $f_2(0,\cdots,0,t_{k+1},\cdots,t_d)$ is positive definite.  If $k=0$ then these are precisely the Hessian conditions mentioned in  Theorem \ref{thr5.4} (3), (4). 
\thr[{\cite[Theorem 9.5.3]{Ma}}] \label{thr5.8} Let $G\subseteq \thuc[X]$ and $f\in \thuc[X]$. Assume
\it 
\item[(1)] $M_G$ is Archimedean;
\item[(2)] $f\geq 0$ on $K_G$;
\item[(3)] each zero $p$ of $f$ in $K_G$ is a non-singular point of $\thuc^{d}$, and there exist $g_1,\cdots g_k\in M_G, $ $1\leq k\leq d$, which are   part of a system of uniformizing parameters at $p$ such that $f$ satisfies BHC with respect to $g_1,\cdots,g_k$ at $p$.
\hit
Then $f \in M_G$.
\ethr
 Note that in this theorem $G$ is an arbitrary  subset of $\thuc[X]$, not necessarily finite.  Using this theorem, we have the following  boundary Hessian criterion for diagonal polynomial matrices.

\pro \label{pro5.9} Let $\hoa{G} \subseteq \hoa{S}_n(\thuc[X])$ and $G\subseteq \thuc[X]$ as in Lemma \ref{lm1.3}.  
Let $\bold{D} = D(d_1,\cdots,d_r)$, $r\leq n$,  be an $n\times n$ diagonal polynomial matrix in $\hoa{S}_n(\thuc[X])$. Assume 
\it  
\item[(1)] $M_G$ is Archimedean;
\item[(2)] $\bold{D}\geq 0$ on $\bold{K}_{\hoa{G}}$;
\item[(3)] each zero $p$ of  each $d_i$ in $\bold{K}_{\hoa{G}}$ is a non-singular point of $\thuc^{d}$, and there exist $g_{i_1},\cdots g_{i_k}\in M_G, $ $1\leq k\leq d$, which are   part of a system of uniformizing parameters at $p$ such that $d_i$ satisfies BHC with respect to $g_{i_1},\cdots,g_{i_k}$ at $p$.
\hit
Then $\bold{D} \in (M_G)^n \subseteq \hoa{M}_{\hoa{G}}$. 
\epro
\pf  The result follows from Theorem \ref{thr5.8}, applying for each $d_i \in \thuc[X]$, and Remark \ref{rem1.2}.
\epf

By a similar argument to the one given  in the proof of  Theorem \ref{thr5.3}, using Proposition \ref{pro5.9}, we obtain the following

\thr \label{thr5.10} Let $\hoa{G}\subseteq \hoa{S}_n(\thuc[X])$ and $G\subseteq \thuc[X]$ as in Lemma \ref{lm1.3}. 
Let $\bold{F}\in \hoa{S}_n(\thuc[X])$. Assume
\it
\item[(1)] $M_G$ is Archimedean;
\item[(2)] $\bold{F}\succcurlyeq \bold{0}$ on $\bold{K}_{\hoa{G}}$;
\item[(3)] for each $\bold{x} \in \thuc^n\tru \{0\}$, each zero $p$ of  the polynomial  $\bold{x}^T\bold{F}\bold{x} $  in $\bold{K}_{\hoa{G}}$ is a non-singular point of $\thuc^{d}$, and there exist $g_{1},\cdots g_{k}\in M_G, $ $1\leq k\leq d$, which are   part of a system of uniformizing parameters at $p$ such that $\bold{x}^T\bold{F}\bold{x} $ satisfies BHC with respect to $g_{1},\cdots,g_{k}$ at $p$.
\hit 
Then 
\it 
\item[(i)] there exists a matrix $\bold{X}_{-}\in \hoa{M}_{n}(\thuc[X])$ such that $\bold{X}_{-}\bold{F}\bold{X}_{-}^T \in (M_G)^n\subseteq \hoa{M}_{\hoa{G}}$;
\item[(ii)] there exists a non-zero polynomial $b\in \thuc[X]$ such that $b^2\bold{F} \in (M_G)^n \subseteq \hoa{M}_{\hoa{G}}$. 
\hit
\ethr

\textbf{Acknowledgements }  The author would like to  thank the anonymous referees for their useful comments and suggestions.  This research is funded by Vietnam National Foundation for Science and Technology Development (NAFOSTED) under the grant number 101.99-2013.24. This work is finished during the author's postdoctoral fellowship at the Vietnam Institute for Advanced Study in Mathematics (VIASM). He thanks VIASM for financial support and hospitality.


\begin{thebibliography}{}
\bibitem  {Ci0} Cimprič, J.: A representation theorem for Archimedean quadratic modules on $*$-rings. Canad. Math. Bull  \textbf{52} (1), 39-52 (2009)
\bibitem  {Ci} Cimprič, J.: Real algebraic geometry for matrices over commutative rings.  J. Algebra \textbf{359}, 
89-103 (2012) 
\bibitem   {CZ} Cimprič, J.,    Zalar, J.:  Moment problems for operator polynomials.   J.  Math.  Anal.  App.  \textbf{401}(1),  307-316 (2013) 
\bibitem  {GR} Gondard, D., Ribenboim, P.: Le 17e probl\`{e}me de Hilbert pour les matrices.   Bull. Sci. Math. (2) \textbf{98} (1),  49–56 (1974) 
\bibitem {JP} Jacobi, T.,   Prestel, A.: Distinguished representations of strictly positive polynomials.   J. reine angew. Math. \textbf{532}, 223-235 (2001)
\bibitem {KSch} Klep, I., Schweighofer, M.:  Pure states, positive matrix polynomials and sums of Hermitian squares.  Indiana
 Uni.  Math. J.  \textbf{59} (3), 857-874 (2010)
\bibitem   {Kr}  Krivine, J.-L.:  Anneaux pr\'{e}odonn\'{e}s.  J. Anal. Math. \textbf{12}, 307–326 (1964)
\bibitem    {Ma} Marshall, M.:  Positive polynomials and sums of squares.   Mathematical Surveys and Monographs vol. \textbf{146},  American Mathematical Society, Providence, RI  (2008) 
\bibitem{Sch1} Scheiderer, C.:  Sums of squares on real algebraic curves. Math. Z. \textbf{245},  725-760 (2003)
\bibitem{Sch2} Scheiderer, C.:  Distinguished representations of non-negative polynomials.   J. Algebra \textbf{289},  558-573 (2005) 
\bibitem {SchH} Scherer,  C. W., Hol, C. W. J.:  Matrix sum-of-squares relaxations for robust semi-definite programs. Math. Program. \textbf{107},  no. 1-2, Ser. B, 189–211 (2006)
 \bibitem  {Schm-1} Schm\"{u}dgen, K.:  Unbounded operator algebras and representation theory. Operator Theory: Advances and Applications, 37. Birkh\"{a}user Verlag, Basel-Boston-Berlin (1990)
\bibitem  {Schm0} Schm\"{u}dgen, K.: A strict Positivstellensatz for the Weyl algebra. Math. Ann. \textbf{331}, 779–794 (2005)
\bibitem  {Schm} Schm\"{u}dgen, K.:  Noncommutative real algebraic geometry - some basic concepts and first ideas. In:  Emerging Applications of Algebraic Geometry, IMA Vol. Math. Appl., vol. \textbf{149}, pp. 325-350. Springer, New York (2009) 
\bibitem   {Schw} Schweighofer, M.:  Global Optimization of polynomials using gradient tentacles and sums of squares.   SIAM J. Optimization \textbf{17}(3), 920-942 (2006) 
\bibitem   {St} Stengle, G.: A Nullstellensatz and a positivstellensatz in semialgebraic geometry. Math. Ann. \textbf{207}, 87–97 (1974) 
\end{thebibliography}


%
%

\end{document}